\documentclass[final]{siamltex}

\newcommand{\beq}{\begin{equation}}
\newcommand{\eeq}{\end{equation}}

\newcommand{\veps}{\varepsilon}

\newcommand{\cx}{{\tilde x}}

\newcommand{\cQ}{{\tilde Q}}
\newcommand{\cR}{{\tilde R}}

\newcommand{\cU}{{\tilde U}}
\newcommand{\hA}{{\widehat A}}
\newcommand{\hQ}{{\widehat Q}}
\newcommand{\hb}{{\widehat b}}

\newcommand{\calA}{{\cal A}}

\newcommand{\beqn}{\[}
\newcommand{\eeqn}{\]}

\newcommand{\bw}{\mbox{\bf w}}

\newcommand{\bs}{\mbox{\bf s}}
\newcommand{\bt}{\mbox{\bf t}}
\newcommand{\bo}{\mbox{\bf 0}}

\newcommand{\boldy}{\mbox{\bf y}}

\title{Stability of Fast Algorithms for\\
	Structured Linear Systems\thanks{%
Appeared as Report TR-CS-97-18, ANU, September 1997.
Preliminary version of a
chapter to appear in
{\em Stability of Fast Methods for Linear Systems with Structure}
(editors, Ali H.~Sayed and Thomas Kailath), SIAM, Philadelphia, 1999.}}

\author{Richard P. Brent\thanks{%
	Copyright \copyright\ 1997, R. P. Brent \hfill rpb177tr}
	}

\begin{document}

\maketitle

\begin{abstract}
We survey the numerical stability of some fast algorithms for solving
systems of linear equations and linear least squares problems with a low
displacement-rank structure. For example, the matrices involved may be
Toeplitz or Hankel.  We consider algorithms which incorporate
pivoting without destroying the structure, %
and describe some recent results
on the stability of these algorithms.
We also compare these results with the corresponding stability
results for the well known algorithms of Schur/Bareiss and Levinson,
and for algorithms based on the semi-normal equations.
\end{abstract}



\pagestyle{myheadings}
\thispagestyle{plain}
\markboth{R. P. BRENT}{NUMERICAL STABILITY OF SOME FAST ALGORITHMS}

\section{Motivation}

\thispagestyle{empty}
The standard direct method for solving dense $n \times n$ systems of
linear equations is Gaussian elimination with partial pivoting.
The usual implementation requires of order $n^3$ %
arithmetic operations. %

In practice, linear systems often arise from some physical system
and have a structure which is a consequence of the physical system.
For example,
time-invariant physical systems often give rise to {\em Toeplitz} systems
of linear equations (\S\ref{subsec:Toeplitz}).
An $n \times n$ Toeplitz matrix is a dense matrix because it generally
has $n^2$ nonzero elements.
However, it is determined by only $O(n)$ parameters (in fact, by the $2n-1$
entries in its first row and column).
Similar examples are Hankel, Cauchy, Toeplitz-plus-Hankel, and Vandermonde
matrices~\cite{Golub89,Olshevsky97a}.

When solving such a structured linear system
it is possible to ignore the structure, and this may have advantages
if standard software is available and $n$ is not too large.
However, if $n$ is large or if many systems have to be solved,
perhaps with real-time constraints (e.g.~in radar and sonar applications),
then it is desirable to take advantage of the structure.
The primary advantage to be gained is that the time to solve a linear
system is reduced by a factor of order $n$ to $O(n^2)$.
Storage requirements may also be reduced by a factor of order $n$,
from $O(n^2)$ to $O(n)$.

Most papers concerned with algorithms for structured linear systems
concentrate on the speed (usually measured in terms of the number of
arithmetic operations required) and ignore questions of numerical
accuracy. However, it is dangerous to use fast algorithms without
considering their numerical properties.  There is no point in obtaining
an answer quickly if it is much less accurate than is justified by
the data.

In this paper we reverse the usual emphasis on speed,
and concentrate on the numerical properties of fast algorithms.
Because there are many classes of structured matrices, and an ever increasing
number of fast algorithms, we can not attempt
to be comprehensive. Our aim is to introduce the reader to the subject,
illustrate some of the main ideas,
and provide pointers to the literature.
\pagebreak[4]

The subject of numerical stability/instability of fast algorithms
is confused for several reasons:
\begin{remunerate}	%
\item Structured matrices are often very ill-conditioned. For example, the
Hilbert matrix~\cite{Hilbert1894},
defined by $a_{i,j} = 1/(i+j-1)$, %
is often used as an example of a poorly conditioned matrix~\cite{Forsythe67}.
This is a typical example of a Hankel matrix. Reversing the order of the rows
gives a Toeplitz matrix.
\item The solution might be less sensitive to structured perturbations
(i.e. perturbations which are physically plausible because they preserve
the structure) than to general (unstructured) perturbations.
The effect of rounding errors in methods which ignore the structure
is generally equivalent to the introduction of unstructured perturbations.
Ideally we should use a method which
only introduces structured perturbations, %
but this property often does not hold, or is difficult to prove,
even for methods which take advantage of the structure to reduce
the number of arithmetic operations.
\item The error bounds which can be proved are usually much weaker than
what is observed on ``real'' or ``random'' examples.  Thus, methods which
are observed to work well in practice can not always be guaranteed,
and it is hard to know if it is just the analysis which is weak or
if the method fails in some rare cases.
\item An algorithm may perform well on special classes of structured matrices,
e.g.~positive definite matrices or matrices with positive reflection
coefficients, but perform poorly or break down on broader classes of
structured matrices.
\end{remunerate}

\subsection{Outline}

Different authors have given different (and sometimes inconsistent)
definitions of {\em stability} and {\em weak stability}.
We shall follow Bunch~\cite{Bunch85,Bunch87}.
For completeness, our definitions are given in~\S\ref{sec:stability}.

The concept of {\em displacement rank},
discussed in \S\ref{sec:classes}, may be used to unify the discussion
of many algorithms for structured
matrices~\cite{Kailath94,Kailath79,Kailath-Sayed}.
It is well known that systems of $n$ linear equations with a low displacement
rank (e.g.\ Toeplitz or Hankel matrices) can be solved in $O(n^2)$ arithmetic
operations. Asymptotically faster algorithms with time bound
$O(n \log^2 n)$ exist~\cite{rpb059},
but are not considered here because their numerical properties are
generally poor and the constant factors hidden in the ``$O$'' notation
are large~\cite{Ammar88}.

For positive definite Toeplitz matrices, the first $O(n^2)$
algorithms were introduced by Kolmogorov~\cite{Kolmogorov41},
Wiener~\cite{Wiener49}
and Levinson~\cite{Levinson47}.
These algorithms are related to recursions of Szeg\"o~\cite{Szego39} for
polynomials orthogonal on the unit circle.  Another class of $O(n^2)$
algorithms, e.g.~the Bareiss algorithm~\cite{Bareiss69}, are related to
Schur's algorithm for finding the continued fraction representation of a
holomorphic function in the unit disk~\cite{Schur17}. This class can be
generalized to cover unsymmetric matrices and more general ``low displacement
rank'' matrices~\cite{Kailath-Sayed}.
In \S\S\ref{sec:structGE}--\ref{sec:fastorthog} we consider the
numerical stability of some of these algorithms.
The GKO-Cauchy and GKO-Toeplitz algorithms are discussed
in~\S\ref{sec:structGE}.
The Schur/Bareiss algorithm for positive definite matrices is
considered in \S\ref{subsec:BareissPD}, and generalized Schur algorithms
are mentioned in \S\ref{subsec:genSchur}.
In \S\ref{sec:fastorthog} we consider fast orthogonal factorization
algorithms and the fast solution of structured least squares problems.

Algorithms for Vandermonde and many other classes of
structured matrices are not considered here~-- we refer to the
extensive survey by Olshevsky~\cite{Olshevsky97a}.
Also, we have omitted any discussion of fast ``lookahead''
algorithms~\cite{Chan92a,Chan92b,Freund93a,Freund93b,Gutknecht92,%
Gutknecht93a,Gutknecht93b,Hansen93,Sweet93}
because, although such algorithms often succeed
in practice, in the worst case they require of order $n^3$ operations.
\pagebreak[4]

Much work has been done on iterative methods for Toeplitz and
related systems. Numerical stability is not a major problem with
iterative methods, but the speed of convergence depends on a good
choice of preconditioner.
We do not consider iterative methods here, except to mention
iterative refinement as a way of improving the accuracy of
weakly stable direct methods.

\subsection{Notation}
\label{subsec:notation}

In the following,
$R$ denotes a structured matrix,
$T$ is a Toeplitz or Toeplitz-type matrix,
$P$ is a permutation matrix,
$L$ is lower triangular,
$U$ is upper triangular, and
$Q$ is orthogonal.
In error bounds, $O_n(\veps)$ means $O(\veps f(n))$,
where $f(n)$ is a polynomial in $n$.
We do not usually try to specify the polynomial $f(n)$ precisely
because it depends on the norm that is used to measure the error
and on many unimportant details of the implementation.

\section{Stability and Weak Stability} %
\label{sec:stability}

In this section we give definitions of stability and weak stability
of algorithms for solving linear systems.

Consider algorithms for solving a nonsingular, $n \times n$ linear
system $Ax = b$.
There are many definitions of %
numerical stability in the literature, for example
\cite{Bjorck87,Bjorck91,rpb144,Bojanczyk91,Bunch85,Cybenko80,%
Golub89,Jankowski77,Miller80,Paige73,Stewart73}.
Definitions~\ref{def:stable} and~\ref{def:weak} below are taken from
Bunch~\cite{Bunch87}. %

\begin{definition}
\label{def:stable}
An algorithm for solving linear equations
is {\em stable} for a class of matrices $\calA$ if for each $A$ in
$\calA$ and for each $b$ the computed solution $\cx$ to $Ax = b$
satisfies $\hA\cx = \hb$, where $\hA$ is close to $A$ and $\hb$ is
close to $b$.
\end{definition}

Definition~\ref{def:stable} says that, for stability,
the {\em computed} solution
has to be the {\em exact} solution of a problem
which is close to the original problem.
This is the classical {\em backward stability} of
Wilkinson~\cite{Wilkinson61,Wilkinson63,Wilkinson65}.
We interpret ``close'' to mean close
in the relative sense in some norm, i.e.
\beqn
	\|\hA - A\|/\|A\| = O_n(\veps),\;
	\|\hb - b\|/\|b\| = O_n(\veps).			%
\eeqn

Note that the matrix $\hA$ is not required to be in the class $\calA$.
For example, $\calA$ might be the class of nonsingular Toeplitz matrices,
but $\hA$ is not required to be a Toeplitz matrix. If we require
$\hA \in \calA$ we get what Bunch~\cite{Bunch87} calls {\em strong stability}.
For a discussion of the difference between stability and strong stability
for Toeplitz algorithms, see~\cite{Gohberg93,Gohberg94,Higham92,Varah92}.

Stability does not imply that the computed
solution $\cx$ is close to the exact solution $x$, unless the problem
is well-conditioned.
Provided $\kappa\veps$ is sufficiently small, stability implies that
\beq
\|\cx - x\|/\|x\| = O_n(\kappa\veps).			\label{eq:relerr}
\eeq
For more precise results,
see~Bunch~\cite{Bunch87} and~Wilkinson~\cite{Wilkinson61}.

As an example, consider the method of Gaussian elimination.
Wilkinson~\cite{Wilkinson61} shows that
\beqn
\|\hA - A\|/\|A\| = O_n(g\veps),			%
\eeqn
where $g = g(n)$ is the ``growth factor''.
$g$ depends on whether partial or complete pivoting is used.
In practice $g$ is usually moderate, even for partial pivoting.
However, a well-known example shows that $g(n) = 2^{n-1}$ is possible
for partial pivoting, and
it has been shown that examples where $g(n)$ grows exponentially
with $n$ may arise in applications, e.g.~for linear systems arising from
boundary value problems.	%
Even for complete pivoting,
it has not been {\em proved} that $g(n)$ is bounded
by a polynomial in $n$.
Wilkinson~\cite{Wilkinson61}	%
showed that $g(n) \le n^{(\log n)/4 + O(1)}$,
and Gould~\cite{Gould91} showed that $g(n) > n$ is possible
for $n > 12$; there is still a large gap between these results.
Thus, to be sure that Gaussian elimination satisfies
Definition~\ref{def:stable}, we must restrict $\calA$
to the class of matrices for which $g$ is $O_n(1)$.
In practice this is not a problem, because $g$ can easily be checked
{\em a~posteriori}~\cite{Wilkinson65}.

Although stability is desirable, it is more than we can prove for many
useful algorithms. Thus, following Bunch~\cite{Bunch87}, we define the
(weaker, but still useful)
property of {\em weak stability}.

\begin{definition}
\label{def:weak}
An algorithm for solving linear equations
is {\em weakly stable} for a class of matrices $\calA$ if for each
well-conditioned $A$ in
$\calA$ and for each $b$ the computed solution $\cx$ to $Ax = b$
is such that $\|\cx-x\|/\|x\|$ is small.
\end{definition}

In Definition~\ref{def:weak}, we take ``small'' to mean $O_n(\veps)$,
and ``well-conditioned'' to mean that $\kappa(A)$ is $O_n(1)$,
i.e.~is bounded by a polynomial in $n$.
From~(\ref{eq:relerr}), stability implies weak stability.

Define the {\em residual} $r = A\cx - b$.
It is well-known~\cite{Wilkinson63} that
\beq
{1 \over \kappa}{\|r\|\over\|b\|} \le
{\|\cx-x\| \over \|x\|} \le \kappa{\|r\|\over\|b\|} .
\eeq
Thus, for well-conditioned $A$, $\|\cx-x\|/\|x\|$ is small if and only if
$\|r\|/\|b\|$ is small. This observation clearly
leads to an alternative definition of weak stability:

\begin{definition}
\label{def:weak2}
An algorithm for solving linear equations
is {\em weakly stable} for a class of matrices $\calA$ if for each
well-conditioned $A$ in
$\calA$ and for each $b$ the computed solution $\cx$ to $Ax = b$
is such that $\|A\cx - b\|/\|b\|$ is small.
\end{definition}

\subsection{Example: orthogonal factorization}

To illustrate the concepts of stability and weak stability,
consider computation of the Cholesky factor $R$ of $A^TA$, where
$A$ is an $m \times n$ matrix of full rank $n$.
A good $O(mn^2)$ algorithm is
to compute the $QR$ factorization
\[	A = QR	\]
of $A$ using Householder or Givens transformations~\cite{Golub89}.
It can be shown~\cite{Wilkinson63}
that the computed matrices $\cQ$, $\cR$ satisfy
\beq
	\hA = \hQ\cR				\label{eq:QRerror}
\eeq
where $\hQ^T\hQ = I$,
$\cQ$ is close to $\hQ$, and $\hA$ is close to $A$.
Thus, the algorithm is stable in the
sense of backward error analysis.
Note that $\|A^TA - \cR^T\cR\|/\|A^TA\|$ is small,
but $\|\cQ - Q\|$ and $\|\cR - R\|/\|R\|$ are not necessarily small.
Bounds on $\|\cQ - Q\|$ and $\|\cR - R\|/\|R\|$
depend on $\kappa$, and
are discussed in~\cite{Golub65,Stewart77,Wilkinson65}.

A different algorithm is to compute (the upper triangular part of) $A^TA$,
and then compute the Cholesky factorization of $A^TA$ by the usual
(stable) algorithm. The computed result $\cR$ is such that
$\cR^T\cR$ is close to $A^TA$. However, this does not
imply the existence of $\hA$ and $\hQ$ such that~(\ref{eq:QRerror}) holds
(with $\hA$ close to $A$ and some $\hQ$ with $\hQ^T\hQ = I$)
unless $A$ is well-conditioned~\cite{Stewart79}. By analogy with
Definition~\ref{def:weak2} above,
we may say that Cholesky factorization of $A^TA$ gives a
{\em weakly stable} algorithm for computing $R$,
because the ``residual'' $A^TA - \cR^T\cR$ is small.

\section{Classes of Structured Matrices}
\label{sec:classes}

Structured matrices $R$ satisfy a
{\em Sylvester equation} which has the form
\beq
\nabla_{\{A_f,A_b\}}(R) = A_f R - R A_b = \Phi\Psi\;, \label{eq:sylv}
\eeq
where $A_f$ and $A_b$ have some simple structure (usually banded, with 3 or fewer
full diagonals), $\Phi$ and $\Psi$ are $n \times \alpha$ and $\alpha \times n$ respectively,
and $\alpha$ is some fixed integer.
The pair of matrices
$(\Phi, \Psi)$ is called
the $\{A_f, A_b\}$-{\em generator} of $R$.

$\alpha$ is called the
$\{A_f, A_b\}$-{\em displacement rank} of $R$.
We are interested in cases where $\alpha$ is small (say at most~4).

For a discussion of the history and some variants of~(\ref{eq:sylv}),
see~\cite{Gohberg95} and the references given there,
in particular~\cite{Heinig84,Kailath79}.

\subsection{Cauchy and Cauchy-type matrices}

Particular choices of $A_f$ and $A_b$ lead to definitions of basic classes of
matrices. Thus, for a Cauchy matrix
$$C(\bt,\bs) = \left[ \frac{1}{t_i - s_j} \right]_{ij}\;,$$
we have
   \beqn A_f = D_t = \mbox{diag}(t_1, t_2, \ldots, t_n)\;, \eeqn
   \beqn A_b = D_s = \mbox{diag}(s_1, s_2, \ldots, s_n)  \eeqn
and
   \beqn \Phi^T = \Psi = [1,1,\ldots,1]\;. \eeqn

As a natural generalization, we can take
$\Phi$ and
$\Psi$ to be any rank-$\alpha$ matrices, with
$A_f, A_b$ as above.  Then a matrix $R$ satisfying the Sylvester
equation~(\ref{eq:sylv}) is said to be a {\em Cauchy-type} matrix.

\subsection{Toeplitz matrices}
\label{subsec:Toeplitz}

For a Toeplitz matrix $T = [t_{ij}] = [a_{i-j}]$, we take
\beqn  A_f = Z_1 = \left[ \begin{array}{ccccc}
                         0 &   0    & \cdots & 0 & 1      \\
                         1 &   0    &        &   & 0      \\
                         0 &   1    &        &   & \vdots \\
                    \vdots &        & \ddots &   & \vdots \\
                         0 & \cdots &    0   & 1 & 0
                         \end{array}  \right],\;\;
A_b = Z_{-1} =  \left[ \begin{array}{ccccc}
                         0 &   0    & \cdots & 0 &  -1    \\
                         1 &   0    &        &   &   0    \\
                         0 &   1    &        &   & \vdots \\
                    \vdots &        & \ddots &   & \vdots \\
                         0 & \cdots &    0   & 1 &   0
                         \end{array}  \right],   \eeqn
\beqn \Phi = \left[ \begin{array}{cccc}
                   1  &    0    &   \cdots    &      0        \\
                  a_0 & a_{1-n}+a_1 & \cdots  & a_{-1}+a_{n-1}
                \end{array} \right]^T,            \eeqn
and
\beqn \Psi =\left[ \begin{array}{cccc}
                a_{n-1}-a_{-1} & \cdots & a_1-a_{1-n} & a_0 \\
                       0       & \cdots &   0 &  1 %
                \end{array} \right]\;.            \eeqn
We can generalize to {\em Toeplitz-type} matrices by taking
$\Phi$ and $\Psi$ to be general rank-$\alpha$ matrices.

\section{Structured Gaussian Elimination}
\label{sec:structGE}

Let an input matrix, $R_1$, have the partitioning
\beqn R_1 = \left[ \begin{array}{cc} d_1 & \bw_1^T \\ \boldy_1 & \tilde{R}_1
       \end{array} \right]\;.
\eeqn
The first step of normal Gaussian elimination is to premultiply $R_1$ by
\beqn
       \left[ \begin{array}{cc} 1 & \bo^T \\ -\boldy_1/d_1 & I
       \end{array} \right]\;,
\eeqn
which reduces it to
\beqn
       \left[ \begin{array}{cc} d_1 & \bw_1^T \\ \bo & R_2
       \end{array} \right]\;,
\eeqn
where
\beqn R_2 = \tilde{R_1} - \boldy_1\bw_1^T/d_1 \eeqn
is the {\em Schur complement}
of $d_1$ in $R_1$.
At this stage, $R_1$ has the factorization
\beqn R_1 = \left[ \begin{array}{cc} 1 & \bo^T \\ \boldy_1/d_1 & I
       \end{array} \right]
       \left[ \begin{array}{cc} d_1 & \bw_1^T \\ \bo & R_2
       \end{array} \right]\;. \eeqn
One can proceed recursively with the Schur complement
$R_2$,
eventually obtaining a factorization $R_1 = LU$.

The key to {\em structured} Gaussian elimination is the fact that
the displacement structure is preserved under Schur complementation, and
that the generators of the Schur complement $R_{k+1}$ can be computed from
the generators of $R_k$ in $O(n)$ %
operations.

Row and/or column interchanges
destroy the structure of matrices such as Toeplitz matrices.
However, if $A_{f}$
is diagonal (which is the case for
Cauchy and Vandermonde	%
type matrices), then {\em the structure is preserved
under row permutations}.

This observation leads to the {\em GKO-Cauchy} algorithm
of Gohberg, Kailath and Olshevsky~\cite{Gohberg95} for fast
factorization of Cauchy-type matrices with partial pivoting, and
many recent variations on the theme by
Boros, Gohberg, Ming Gu, Heinig, Kailath, Olshevsky, M.~Stewart, {\em et~al}:
see \cite{Boros95,Gohberg95,Gu95a,Heinig94,Olshevsky97a,Stewart97a}.

\subsection{The GKO-Toeplitz algorithm}

Heinig~\cite{Heinig94} showed that, if $T$ is a Toeplitz-type matrix, then
\beqn
R = FTD^{-1}F^\ast
\eeqn
is a Cauchy-type matrix,
where
\beqn F = \frac{1}{\sqrt n} [ e^{2\pi i(k-1)(j-1)/n} ]_{1\leq k,j\leq n} \eeqn
is the Discrete Fourier Transform matrix,
\beqn
D = \diag(1,e^{\pi i/n},\ldots,e^{\pi i(n-1)/n}),
\eeqn
and the generators of $T$ and $R$ are simply related. %

The transformation $T \leftrightarrow R$ is perfectly stable
because $F$ and $D$ are unitary.
Note that $R$ is (in general) complex even if $T$ is real.
This increases the constant factors in the time bounds, because
complex arithmetic is required.

Heinig's observation was exploited by
Gohberg, Kailath and Olshevsky~\cite{Gohberg95}:
$R$ can
be factorized as $R = P^TLU$ using GKO-Cauchy.
Thus, from the factorization
\beqn
T = F^\ast P^T L U F D  \;,
\eeqn
a linear system involving $T$ can be solved in $O(n^2)$ operations.
The full procedure of conversion to Cauchy form, factorization,
and solution requires $O(n^2)$ (complex) operations.

Other structured matrices, such as
Hankel,
Toeplitz-plus-Hankel,
Vandermonde,
Chebyshev-Vandermonde, etc, can be converted
to Cauchy-type matrices in a similar way.

\subsection{Error Analysis}

Because GKO-Cauchy and GKO-Toeplitz involve partial pivoting,
we might guess that their stability would be similar to that of
Gaussian elimination with partial pivoting.
Unfortunately, there is a flaw in this reasoning. During
GKO-Cauchy the {\em generators} have to be transformed, and the
partial pivoting does not ensure that the transformed generators
are small.

Sweet and Brent~\cite{rpb157}
show that significant generator growth can occur if all
the elements of $\Phi\Psi$ are
small compared to those of $|\Phi||\Psi|$.
This can not happen for ordinary Cauchy matrices because
$\Phi^{(k)}$ and $\Psi^{(k)}$ have only one column and one
row respectively. However, it can happen for higher displacement-rank
Cauchy-type matrices, even if the original matrix
is well-conditioned.

\subsection{The Toeplitz Case}

In the Toeplitz case there is an extra constraint on the selection of
$\Phi$ and $\Psi$, but it is still possible to give examples where
the normalized solution error grows like $\kappa^2$ and the
normalized residual grows like $\kappa$,
where $\kappa$ is the condition number of the Toeplitz matrix.
Thus, the GKO-Toeplitz
algorithm is (at best) weakly stable.

It is easy to think of modified algorithms which avoid the examples
given by Sweet and Brent, but it is difficult to prove that they are stable
in all cases. Stability depends on the worst case, which may be rare and
hard to find by random sampling.

The problem with the original GKO algorithm is growth in the generators.
Ming Gu suggested exploiting the fact that the generators are not unique.
Recall the
{Sylvester equation} (\ref{eq:sylv}).
Clearly we can replace $\Phi$ by $\Phi M$ and
$\Psi$ by $M^{-1} \Psi$, where $M$ is any invertible $\alpha \times \alpha$
matrix, because this does not change the product $\Phi\Psi$.
Similarly at later stages of the GKO algorithm.

Ming Gu~\cite{Gu95a}
proposes taking $M$ to orthogonalize the columns of $\Phi$
(that is, at each stage perform an orthogonal factorization of the
generators).
M.~Stewart~\cite{Stewart97a}
proposes a (cheaper) LU factorization of the generators.
In both cases, clever pivoting schemes give error bounds analogous to
those for Gaussian elimination with partial pivoting.

\subsection{Gu and Stewart's error bounds}

The error bounds obtained by Ming Gu and M.~Stewart involve
a factor $K^n$, where $K$ depends on the ratio of
the largest to smallest modulus elements in the Cauchy matrix
\beqn \left[ \frac{1}{t_i - s_j} \right]_{ij}\;.\eeqn
Although this is unsatisfactory, it is similar to the factor $2^{n-1}$
in the error bound for Gaussian elimination with partial pivoting.
In practice, the latter factor is extremely pessimistic, which explains why
Gaussian elimination with partial pivoting is
popular~\cite{Higham96}. Perhaps the bounds of Ming Gu and
M.~Stewart are similarly pessimistic, although practical experience
is not yet extensive enough to be confident of this.
M.~Stewart~\cite{Stewart97a}
gives some interesting numerical results which indicate that his scheme
works well, but more numerical experience is necessary
before a definite conclusion can be reached.

\subsection{A general strategy for guaranteed results}

It often happens that there is a choice of
\begin{remunerate}
\item A fast algorithm which usually gives
an accurate result, but occasionally fails (or at least can not be proved
to succeed every time), or
\item An algorithm which is guaranteed to be stable but is slow.
\end{remunerate}
In such cases, a good strategy is to use the fast algorithm, but then
check the normalized residual $||A\cx - b||/(||A||\cdot||\cx||$,
where $\cx$ is the computed solution of the system $Ax = b$.
If the residual is sufficiently small we can accept $\cx$ as a reasonably
accurate solution.  In the rare cases that the residual is not sufficiently
small we can can use the slow but stable algorithm (alternatively, if
the residual is not too large, one or two iterations of iterative
refinement may be sufficient and faster~\cite{Jankowski77,Wilkinson65}).

An example of this general strategy is the solution of a Toeplitz system
by Ming Gu or M.~Stewart's modification of the GKO algorithm.
We can use the $O(n^2)$ algorithm, check the residual,
and resort to iterative refinement or a stable $O(n^3)$ algorithm
in the (rare) cases that it is necessary.
Computing the residual takes only $O(n\log n)$ arithmetic operations.

\section{Positive Definite Structured Matrices}
\label{sec:PosDef}

An important class of algorithms, typified by the algorithm of
Bareiss~\cite{Bareiss69},
find an $LU$ factorization of a Toeplitz matrix $T$,
and (in the symmetric case)
are related to the classical algorithm of
Schur~\cite{Burg75,Gohberg86,Schur17}.

It is interesting to consider the numerical properties of these algorithms
and compare with the numerical properties of the Levinson algorithm
(which essentially finds an $LU$ factorization of $T^{-1}$).

\subsection{The Bareiss algorithm for positive definite matrices}
\label{subsec:BareissPD}

Bojanczyk, Brent, de Hoog and
Sweet (abbreviated BBHS) have shown in~\cite{rpb144,Sweet82}
that the numerical properties of the Bareiss
algorithm are similar to
those of Gaussian elimination ({\em without} pivoting).
Thus, the algorithm is stable for positive definite symmetric
Toeplitz matrices.

The Levinson algorithm can be shown to be weakly stable for bounded $n$,
and numerical results by Varah~\cite{Varah93},
BBHS and others suggest that this is all that
we can expect. Thus, the Bareiss algorithm is (generally) better numerically
than the Levinson algorithm.

Cybenko~\cite{Cybenko80} showed that if certain quantities called
``reflection coefficients'' are positive %
then the Levinson-Durbin
algorithm for solving the
Yule-Walker equations %
(a positive-definite system
with special right-hand side) is stable.
Unfortunately, Cybenko's result is not usually applicable, because
most positive-definite Toeplitz matrices do not usually
satisfy the restrictive condition on the reflection coefficients.
Some relevant numerical examples are given in~\cite{rpb144}.

\subsection{The generalized Schur algorithm}
\label{subsec:genSchur}

The Schur algorithm can be generalized to factor a large variety of
structured matrices~\cite{Kailath94,Kailath-Sayed}.
For example, suitably generalized Schur algorithms
apply to block Toeplitz matrices, Toeplitz block matrices,
and to matrices of the form $T^T T$, where $T$ is rectangular Toeplitz.

It is natural to ask if the stability results of BBHS (which are for the
classical Schur/Bareiss algorithm) extend to the generalized Schur algorithm.
This was considered by
Chandrasekharan and Sayed~\cite{Chandra96a} and
M.~Stewart and Van Dooren~\cite{Stewart97b}.
The displacement structure considered in~\cite{Chandra96a} is
more general than that considered in~\cite{Stewart97b},
which in turn is more general than that considered in~\cite{rpb144,Sweet82}.
Each increase in generality complicates the analysis and seems to
make the stability result more dependent on details of the
method of implementing the hyperbolic rotations which occur in the algorithms.
The interested reader is referred to~\cite{Chandra96a,Stewart97b} for
details.

The overall conclusion is that the generalized Schur algorithm is stable for
positive definite symmetric (or Hermitian) matrices,
{\em provided} the hyperbolic transformations
in the algorithm are implemented correctly.

\section{Fast Orthogonal Factorization}
\label{sec:fastorthog}

In an attempt to achieve stability without pivoting,
and to solve $m \times n$ least
squares problems ($m \ge n$), it is natural
to consider algorithms for computing an orthogonal factorization
\beq
	T = QU					\label{eq:AQR}
\eeq
of an $m \times n$ Toeplitz matrix $T$. We assume that $T$ has full rank~$n$.
For simplicity, in the time bounds we assume $m = O(n)$ to avoid
functions of both $m$ and $n$.

The first $O(n^2)$ (more precisely, $O(mn)$) algorithm
for computing the factorization~(\ref{eq:AQR})
was introduced by
Sweet~\cite{Sweet82}. Unfortunately, Sweet's algorithm
is unstable: see Luk and Qiao~\cite{Luk87}.

Other $O(n^2)$ algorithms for computing the matrices $Q$ and $U$ or $U^{-1}$
were given by
Bojanczyk, Brent and de Hoog (abbreviated BBH)~\cite{rpb092},
Chun {\em et~al}~\cite{Chun87}, %
Cybenko~\cite{Cybenko87},
and Qiao~\cite{Qiao88},
but none of them has been shown
to be stable, %
and in several cases examples show that they are unstable.

It may be surprising that fast algorithms for computing an orthogonal
factorization~(\ref{eq:AQR}) are unstable.
The classical $O(n^3)$ algorithms
are stable because they form $Q$ %
as a product of elementary orthogonal matrices (usually Givens or
Householder matrices~\cite{Golub65,Golub89,Wilkinson65}).
Unlike the classical algorithms,
the $O(n^2)$ algorithms do not form $Q$
in a numerically stable manner
as a product of matrices which are (close to) orthogonal.
This observation explains both their speed and their instability!

For example, the algorithms of BBH~\cite{rpb092}
and Chun {\em et~al}~\cite{Chun87} depend on Cholesky downdating,
and numerical experiments show that they do not give a $Q$ which is
close to orthogonal. This is not too surprising, because Cholesky
downdating is known to be a sensitive numerical
problem~\cite{rpb095,Stewart79}.

\subsection{Use of the semi-normal equations}

It can be shown that,
provided the Cholesky downdates are implemented
in a certain way
(analogous to the condition for the stability of the
generalized Schur algorithm),
the BBH algorithm computes $U$ in a weakly stable manner~\cite{rpb143}.
In fact, the computed upper triangular matrix $\cU$
is about as good as can be obtained by performing a Cholesky
factorization of $T^TT$, so
\beqn \|T^TT - \cU^T\cU\|/\|T^TT\| = O_m(\veps)\;.\eeqn
Thus, by solving
\beqn \cU^T\cU x = T^Tb\eeqn
(the so-called {\em semi-normal} equations)
we have a {\em weakly stable} algorithm for the solution of general Toeplitz
systems $Tx = b$ in $O(n^2)$ operations.
The solution can be improved by iterative refinement if desired~\cite{Golub66}.
The computation of $Q$ is avoided, and the algorithm
is applicable to full-rank Toeplitz least squares problems.
The disadvantage of this method is that, by implicitly forming $T^T T$,
the condition of the problem is effectively squared.
If the condition number $\kappa = \kappa(T)$ is in the range
\beqn
{1 \over \sqrt{\veps}} \le \kappa \le {1 \over \veps}
\eeqn
then it will usually be impossible to get any significant figures in the
result (iterative refinement may fail to converge) without reverting to
a slow but stable orthogonal factorization algorithm.
One remedy is to use
double-precision arithmetic, i.e.\ replace $\veps$ by $\veps^2$,
but this may be difficult if $\veps$ already corresponds to the
maximum precision implemented by hardware.

Another way of computing the upper triangular matrix $U$, but
not the orthogonal matrix $Q$, in~(\ref{eq:AQR}),
is to apply the generalized Schur algorithm to $T^T T$.
This method also squares the condition number.

\subsection{Computing $Q$ stably}

It seems difficult to give a satisfactory $O(n^2)$ algorithm for
the computation of $Q$ in the factorization~(\ref{eq:AQR}).
Using a modification of the embedding approach pioneered
by Chun and Kailath~\cite{Chun89,Kailath94},
Chandrasekharan and Sayed~\cite{Chandra96b} %
give a stable algorithm to compute a factorization
\beq
T = LQU        						\label{eq:LQU}
\eeq
where $L$ is lower triangular. %
(Of course, the factorization~(\ref{eq:LQU}) is not unique.)
The algorithm of~\cite{Chandra96b} can be used to solve linear equations,
but not least squares problems
(because $T$ has to be square, and in any case
the matrix $Q$ in~(\ref{eq:LQU})
is different from the matrix $Q$ in~(\ref{eq:AQR})).
Because the algorithm involves embedding the $n \times n$ matrix $T$ in
a $2n \times 2n$ matrix
\beqn \left[ \begin{array}{cc} T^T T & T^T \\ T & 0
       \end{array} \right]\;,
\eeqn
the constant factors in the operation count
are large:
$59n^2 + O(n \log n)$,
which should be compared to $8n^2 + O(n \log n)$ for BBH and the semi-normal
equations. (These operation counts apply for $m = n$: see~\cite{rpb143} for
operation counts of various algorithms when $m \ge n$.)
Thus, although the embedding approach is elegant and leads to interesting
(and in some cases stable) $O(n^2)$ algorithms, a penalty is the significant
increase in the constant factors. This is analogous to the penalty paid by
the GKO algorithm because of the introduction of complex arithmetic.

\section{Conclusions}

Although this survey has barely scratched the surface,
we hope that the reader who has
come this far is convinced that questions of numerical stability are
amongst the most interesting, difficult, and useful questions that we
can ask about fast algorithms for structured linear systems.
It is not too hard to invent a new fast algorithm, but to find a new
stable algorithm is more difficult, and to {\em prove} its stability
or weak stability is a real challenge!
\pagebreak[4]

\subsection*{Acknowledgements}

A preliminary version of this review appeared in~\cite{rpb173}.
Thanks to
Greg Ammar,
Adam Bojanczyk,
James Bunch,
Shiv Chandrasekharan,
George Cybenko,
Paul Van Dooren,
Lars Eld\'en,
Roland Freund,
Andreas Griewank,
Ming Gu,
Martin Gutknecht,
Georg Heinig,
Frank de Hoog,
Franklin Luk,
Vadim Olshevsky,
Haesun Park,
Ali Sayed,
Michael Stewart,
Douglas Sweet
and James Varah
for their assistance.

\end{document}